\documentclass[a4paper,12pt]{article}
\usepackage[all]{xy}
\usepackage{amsmath, amsthm}
\usepackage{amssymb,amsfonts}
\usepackage{array}
\usepackage{hyperref}

\title{Eulerian of the Zero Divisor graph $\Gamma[\mathbb {Z}_n]$}
\author{B.Surendranath Reddy, Rupali.S.Jain and N.Laxmikanth\\
\textit{surendra.phd@gmail.com},{rupalisjain@gmail.com} and\\
{laxmikanth.nandala@gmail.com}\\
Swami Ramanand Teerth Marathwada University, Nanded-431606, India.}{}
\date{}
\begin{document}
\theoremstyle{definition}
\newtheorem{definition}{Definition}[section]
\newtheorem{example}[definition]{Example}
\newtheorem{remark}[definition]{Remark}
\newtheorem{observation}[definition]{Observation}
\theoremstyle{plain}
\newtheorem{theorem}[definition]{Theorem}
\newtheorem{lemma}[definition]{Lemma}
\newtheorem{proposition}[definition]{Proposition}
\newtheorem{corollary}[definition]{Corollary}
\newtheorem{AMS}[definition]{AMS}
\newtheorem{keyword}[definition]{keyword}
\maketitle
\section*{Abstract}
The Zero divisor Graph of a commutative ring $R$, denoted by $\Gamma[R]$, is a graph whose vertices are non-zero zero divisors of $R$ and two vertices are adjacent if their product is zero. We consider the  zero divisor graph $\Gamma[\mathbb{Z}_n]$, for any natural number $n$ and find out which graphs are Eulerian graphs.\\
Keywords:-Zero divisor graph, Euler tour, Euler graph.\\
MSC:- 05C12,05C25,05C50.
\section{Introduction}
The concept of the Zero divisor graph of a ring $R$ was first introduced by I.Beck\cite{3} in 1988 and later on Anderson and Livingston\cite{2}, Akbari and Mohammadian\cite{1}  continued the study of zero divisor graph by considering only the non-zero zero divisors. The concepts of the Euler graph  found in\cite{4}. In this paper we introduce the concepts of  the Euler graph  to the  zero divisor graph $\Gamma[{\mathbb{Z}_n}]$ and identify which zero divisors graphs are Eulerian.  \\
 In this article,  section 2, is about the preliminaries and notations related to zero divisor graph of a commutative ring $R$, in section 3, we derive the  Euler graphs of a zero divisor graph  $\Gamma[{\mathbb{Z}_{p^n}}]$,  and in section 4, we discuss about  Euler graphs of $\Gamma[{\mathbb{Z}_n}]$ for any natural number $n $.
\section{Preliminaries and Notations}
\begin{definition}{\textbf{Zero divisor Graph}}{\cite{1,2}}\\
Let R be a commutative ring with unity and $ Z[R]$ be the set of its zero divisors. Then the zero divisor graph of R denoted by  $\Gamma[R]$, is the graph(undirected) with vertex set  $Z^*[R]= Z[R]-\{\mathbf{0}\}$, the non-zero zero divisors of $R$, such that two vertices $v,w \in Z^*[R]$ are adjacent if $vw=0$.
\end{definition}
\begin{definition}{\textbf{Euler tour}}{\cite{4}}\\
An Euler tour of a graph  $G$  is a tour which includes each edge of the graph  $G$ exactly once.
\end{definition}
\begin{definition}{\textbf{Euler graph}}{\cite{4}}\\
 A graph $G$ is called Euler graph or Eulerain if it has an Euler tour.\\
 \begin{theorem}{\cite{4}}
 A connected graph is Euler iff the degree of every vertex is even.
 \end{theorem}
\end{definition}
\section{Eulerian of The zero divisor graph $\Gamma[\mathbb{Z}_{p^n}]$}
In this section, we discuss the Eulerian of the zero divisor graph  $\Gamma[\mathbb{Z}_{p^n}]$ where p is a prime number.\\ To start with, we consider the zero divisor graph $\Gamma[\mathbb{Z}_n]$ for $n=p^2$.
\begin{theorem}
The zero divisor graph $\Gamma[\mathbb{Z}_{p^2}]$ is  a Euler graph if and only if $p>2$.
 \begin{proof}
Consider  zero divisor graph $\Gamma[\mathbb{Z}_{p^2}]$.\\ The vertex set is 
 $A =\{kp\,|\,k= 1,2,3,....,p-1\,  \}$ and so $|A|=(p-1)$.\\
 As product of any two vertices is zero, they are adjacent and so the corresponding graph is a complete graph on (p-1) vertices that is, $\Gamma[\mathbb{Z}_{p^2}]=K_{p-1}$.\\
As the graph is complete,the degree of each and every vertex of it is (p-1).\\
 If $p>2$ then every prime greater than 2 is odd and hence the degree of each vertex is even.Thus $\Gamma[\mathbb{Z}_{p^2}]$ is Eulerian.\\
For $ p=2$ then the corresponding graph has no Euler path as it consists of only one vertex,thus $\Gamma[\mathbb{Z}_{4}]$ is not Eulerian.
\end{proof}
\end{theorem}
\begin{theorem}
The zero divisor graph $\Gamma[\mathbb{Z}_{p^3}]$ is not an Euler graph, for any prime $p$.
\begin{proof}
Consider the zero divisor graph $\Gamma[\mathbb{Z}_{p^3}]$.\\
Here, we divide the elements(vertices) of $\Gamma[\mathbb{Z}_{p^3}]$ into two disjoint sets namely multiples of $p$ and the multiples of $p^2$ which are given  by
\begin{align*}
 A &=\{kp\,|\,k= 1,2,3,....,p^2-1\, \text{and}\, k\nmid p \}  \\
 B &=\{lp^2\,|\,l= 1,2,3,....,p-1 \} 
\end{align*}
with cardinality  $|A|=p(p-1)$ and $|B|=(p-1)$.\\
As every element of $A$ is adjacent only with the elements of $B$, the degree of each and every vertex of $A$ is $(p-1)$ which is even.\\ 
 Also every element of $B$  is adjacent with itself and with every element of $A$.\\
 Therefore the degree of each and every vertex of $B$ is given by $|A|+|B|-1 $ that is $(p^2-2)$ which is odd.\\  
Hence  $\Gamma[\mathbb{Z}_{p^3}]$ is not Eulerian.\\
If $p=2$, then degree of each vertex of $A$ is $p-1$ which is odd.
Therefore the zero divisor graph $\Gamma[\mathbb{Z}_{p^3}]$ is not  an Euler graph.\\
 \end{proof}
\end{theorem}
With similar arguments, we prove the more general case in the following theorem.
\begin{theorem}
The zero divisor graph $\Gamma[\mathbb{Z}_{p^n}]$ is not Eulerian, for any prime $p$.
 \begin{proof}
We divide the elements(vertices) of $\Gamma[\mathbb{Z}_{p^n}]$ into $n-1$ disjoint sets namely multiples of $p$, multiples of $p^2$... multiples of $p^{n-1}$, given by
\begin{align*}
 A_{1} &=\{k_{1}p\,|\,k_{1}= 1,2,3,....,p^{n-1}-1\, \text{and}\, k_{1}\nmid p \}  \\
 A_{2} &=\{k_{2}p^2\,|\,k_{2}= 1,2,3,....,p^{n-2}-1\, \text{and}\, k_{2}\nmid p^2 \}\\ 
A_{i} &=\{k_{i}p^i\,|\,k_{i}= 1,2,3,....,p^{n-i}-1\, \text{and}\, k_{i}\nmid p^i \} 
\end{align*}
with cardinality  $|A_{i}|=(p^{n-i}-p^{n-i-1})$ ,for $i=1,2,......n-1$.\\
Also the smallest set is $A_{n-1}$ of order $p-1$.\\
Now the degree of an element $v_{i}$ in $A_{i}$ is $p^{i}-2$ which is odd $\forall$ $i=[\frac{n}{2}] $ a greatest integer part function,since the elements of $A_{i}$ are adjacent with itself and also with $A_{j}$ for  $j\geq[\frac{n}{2}] $ .\\
We can make a similar argument for all other sets i.e., every element of $A_{i}$  is adjacent with every element of $A_{n-j}$ where $j\leq i$, therefore the degree of every vertex of  $A_{i}$ is      $\sum_{j=1}^{i}(p^{n-j}-p^{n-j-1})-1=(p^{n-1}-2)$ which is odd.\\
Hence the zero divisor graph $\Gamma[\mathbb{Z}_{p^n}]$ is not Eulerian.\\
If $ p=2$, then the degree of an element $v_{1}$ in $A_{1}$ is $p-1$, which is odd.
Therefore, for any prime, the zero divisor graph $\Gamma[\mathbb{Z}_{p^n}]$ is not  Eulerian.
 \end{proof}
\end{theorem}
\section{Eulerian of the zero divisor graph $\Gamma[\mathbb{Z}_n]$}
In this section we discuss the Eulerian of the zero divisor graph $\Gamma[{\mathbb{Z}_n}]$  where $ n=p_{1}^{\alpha_{1}}p_{2}^{\alpha_{2}}.....p_{k}^{\alpha_{k}}$.\\ To start with, we consider $n=pq$.
\begin{theorem}
The zero divisor graph $\Gamma[\mathbb{Z}_{pq}]$ is  a Euler graph iff $p$ and $q2$ are odd.
 \begin{proof}
Consider the zero divisor graph $\Gamma[\mathbb{Z}_{pq}]$.\\
clearly $\Gamma[\mathbb{Z}_{pq}]$ is a complete bipartite graph, the vertex sets are given by
\begin{align*}
 A &=\{kp\,|\,k= 1,2,3,....,p-1\, \text{and}\, k\nmid q \}  \\
 B &=\{lq\,|\,l= 1,2,3,....,q-1 \, \text{and}\, k\nmid p \}  
\end{align*}
with cardinality  $|A|=(p-1)$ and $|B|=(q-1)$.\\
If  $p$ and $q2$ are odd, then $(p-1)$ and $(q-1)$ are even implies the degree of every vertex of the graph is even and thus the respective graph is an Euler graph.\\
If  $p\, \text{or}\, q=2$, then clearly the graph is not Eulerian as the degree of the atleast one vertex is odd.
 \end{proof}
\end{theorem}
\begin{theorem}
The zero divisor graph $\Gamma[\mathbb{Z}_{p^{\alpha}q^{\beta} }]$ is not Eulerian for all $ \alpha ,\beta\neq 1$ .
 \begin{proof}
Consider the zero divisor graph $\Gamma[\mathbb{Z}_{p^{\alpha}q^{\beta} }]$.\\
Here, we divide the vertices of $\Gamma[\mathbb{Z}_{p^{\alpha}q^{\beta} }]$   into  disjoint sets namely multiples of $ p^{i}$, multiples of $q^{j}$ and multiples of $p^{i}q^{j}$ given by
\begin{align*}
 A_{p^{i}} &=\{r_{i}p^{i}\,|\,r_{i}= 1,2,3,....,p^{i}-1\, \text{and}\, r_{i}\nmid p^{i} \}\\
 A_{q^{j}} &=\{s_{j}q^{j}\,|\,s_{j}= 1,2,3,....,q^{j}-1\, \text{and}\, s_{j}\nmid q^{j} \} \\ 
A_{p^{i}q^{j}} &=\{t_{ij}p^iq^{j}\,|\,t_{ij}= 1,2,3,....,p^{i}q^{j}-1\, \text{and}\, t_{ij}\nmid p^i\, \text{and}\, t_{ij}\nmid q^j \}.
\end{align*}
Then the order of the sets are  $|A_{p^{i}} |=(p^{i}-1)$, $|A_{q^{j}} |=(q^{j}-1)$ and $|A_{p^{i}q^{j}}|=(p^{i}-1)(q^{j}-1)$.\\ 
Assume that both $p$ and $q$ are odd primes.\\
Since every element of the set $ A_{p^{i}}$ is adjacent with the elements of $A_{p^{j}}$, the degree of each and every vertex of the set$ A_{p^{i}}$  is $(q^{j}-1)$.\\
Similarly the degree of each and every vertex of the set $ A_{p^{j}}$  is $(p^{i}-1)$ and the degree of each and every vertex of the set $ A_{p^{i}q^{j}}$  is $|A_{p^{i}} |+A_{p^{j}}|+|A_{p^{i}q^{j}}|-1 = (p^{i}-1)+(q^{j}-1)+(p^{i}-1)(q^{j}-1) = (p^{i}q^{j}-2) $,  which is odd.\\
Thus the degree of the vertices of the corresponding sets is odd.\\
Hence the zero divisor graph $\Gamma[\mathbb{Z}_{p^{\alpha}q^{\beta}}]$ is not Eulerian.\\
If one of  $p\, \text{or}\, q=2$, then also the graph is not Eulerian as the degree of the atleast one vertex is $(p^i-1)$ or $(q^j-1)$ which is odd. Also  the  degree of the elements of $A_{p}$ is $(p-1)$ as the elements of these set are adjacent only with the elements of the set  $A_{p^{\alpha-1}q^{\beta}}$.
\end{proof}
\end{theorem}
\begin{theorem}
 The zero divisor graph  $\Gamma[\mathbb{Z}_{n }]$ where $n =p_{1}^{\alpha_{1}}p_{2}^{\alpha_{2}}.....p_{k}^{\alpha_{k}}$ is not Eulerian for $ \alpha_{i}\geq 2$ where          $ i = 1,2,.........,k$ .
 \begin{proof}
Consider a zero divisor graph  $\Gamma[\mathbb{Z}_{n}]$  where $n =p_{1}^{\alpha_{1}}p_{2}^{\alpha_{2}}.....p_{k}^{\alpha_{k}}$  .\\
Here,we divide the elements(vertices) of $\Gamma[\mathbb{Z}_{n}]$  into  the corresponding  disjoint sets of product of all possible powers of given primes like set of powers of $p_j^i$, set of product of powers of $p_i^rp_j^s$ and so on.\\
Among these sets, we consider the sets of the form\\
$A_{i} =\{ m( p_{1}^{\alpha_{1}}p_{2}^{\alpha_{2}}.....p_{i}^{\alpha_{i-1}}p_{i}^{\alpha_{i+1}}..........p_{k}^{\alpha_{k}})\,\}$ with  $|A_{i}|= {p_i^{{\alpha}_i}}-1$.\\
Now consider the set  $ A_{{p^{\alpha}_i}_i} =\{ t p^{\alpha_{i}}_{i}\,|\, t\nmid {p^{\alpha_{i}}_i}\}$.\\
Assume that all the primes are odd.\\
Since the elements of $ A_{p^{i}_j}$ are adjacent only with the vertices of $ A_j $, the degree of each and every vertex of the set $ A_{p^{i}_j}$  is ${p_i^{{\alpha}_i}}-1 $ which is odd. 
Hence the zero divisor graph  $\Gamma[\mathbb{Z}_{n }]$ where $n =p_{1}^{\alpha_{1}}p_{2}^{\alpha_{2}}.....p_{k}^{\alpha_{k}}$ is not Eulerian for $ \alpha_{i} \geq 2$ where     $ i = 1,2,.........,k$.\\
If one of $ p_{i}=2$, then also the graph is not Eulerian as the degree of the atleast one vertex is ${p_i^{{\alpha}_i}}-1$, which is odd.
\end{proof}
\end{theorem}
We conclude that  the zero divisor graph  $\Gamma[\mathbb{Z}_{n}]$ is Eulerain if and only if either $n=p^2$ or $n=pq$.


\begin{thebibliography}{10}
\bibitem[1]{1} S.Akbari, A.Mohammadian, On the Zero divisor graph of a commutative rings, J.Algebra.2004;274:847-855.
\bibitem[2]{2} D.F.Anderson and P.S.Livingston, The Zero divisor graph of Commutative ring, J.Algebra 217(1999),no.2,434-447.
\bibitem[3]{3} I.Beck,  Coloring of Commutative rings, J.Algebra 116(1988), no.1.208-226.
\bibitem[4]{4} Eulerain graphs and Related Topics part 1,by Herbert FLEISCHNER.
\bibitem[5]{5}  R. C. READ, Euler graphs on labelled nodes, Canad. J. Math., 14 (1962), pp. 482-486.
\bibitem[6]{6}  R. W. ROBINSON, Enumeration of Euler graphs, Proof Techniques in Graph Theory, F. Harary,
ed., Academic Press, N.Y., 1969, pp. 147-153.

\end{thebibliography}
\end{document}